\def\versiondate{7 Aug. 1999}

\input math.macros
\input Ref.macros

\newif\ifsubmit
%\submittrue

%\proofmodetrue
%\leftsectionheadtrue
\checkdefinedreferencetrue
%\continuousnumberingtrue
\continuousfigurenumberingtrue
\theoremcountingtrue
\sectionnumberstrue
%\figuresectionnumberstrue
\forwardreferencetrue
%\lefteqnumberstrue
%\tocgenerationtrue
\initialeqmacro

\def\supp{{\rm supp}}
\def\mat(#1,#2,#3,#4){\pmatrix{#1&#2\cr #3&#4\cr}}
\def\vec(#1,#2){\pmatrix{#1\cr#2\cr}}
\def\L{{\cal L}} %% law

\def\u{{\bf u}}
\def\M{{\bf M}}

\def\PP{{\Bbb P}}
\def\tensormat(#1,#2,#3,#4,#5,#6){\pmatrix{#1&#2&#2&#4\cr #1&#3&#2&#5\cr
     #1&#2&#3&#5\cr #1&#3&#3&#6\cr}}

\ifsubmit
\pageno=0
\title{Singularity of Some Random Continued Fractions}
\vfil
Russell Lyons\footnote{*}{
Department of Mathematics,
Indiana University,
Bloomington, IN 47405-5701
\qquad\par
{\tt rdlyons@indiana.edu}
}
\vfil
Running head: Random Continued Fractions
\vfil
Telephone number for Russell Lyons:
812-855-1645; fax 812-855-0046.
\vfil\eject
\fi

\def\firstheader{\eightpoint\ss\underbar{\raise2pt\line 
    {To appear in {\it J. Theoret. Prob.}\hfil Version of \versiondate}}}

\vglue20pt

\beginniceheadline

\title{Singularity of Some Random Continued Fractions}

\author{Russell Lyons}

\abstract{We prove singularity of some distributions of random continued
fractions that correspond to iterated function systems with overlap and a
parabolic point. These arose while studying the conductance of
Galton-Watson trees.
}

\bottom{Primary 60G30.  %Continuity and singularity of induced measures
 Secondary 28A78.} %Hausdorff measures
{Conductance, Galton-Watson trees, iterated function system, linear
fractional map.}
{Research partially supported by the Universit\'e de Lyon and NSF grant
DMS-9802663.}

\bsection {Introduction}{s.intro}

We are interested in the distributions of certain random continued
fractions. They arose while studying the following question about
Galton-Watson trees: Suppose that $p_k \in \CO{0, 1}$ for $k \ge 1$ with
$\sum_k p_k = 1$. Let $T$ be the random genealogical tree resulting from
the associated Galton-Watson process beginning with one individual, $\rho$,
where each individual has $k$ children with probability $p_k$.  Adjoin a
parent $\Delta$ to $\rho$. Then the probability $\gamma$ that simple random
walk starting at $\Delta$ will return to $\Delta$ is equal to the effective
conductance of $T \cup \{\Delta\}$ from $\Delta$ to infinity when each edge
has unit conductance (see, e.g., Doyle and Snell (1984) or Lyons and Peres
(1998)).  This quantity $\gamma$ enters, for example, in calculations of
the Hausdorff dimension of harmonic measure on the boundary of $T$ (see
Lyons, Pemantle and Peres (1995)).  Thus, it is of interest to calculate
the distribution of $\gamma$. If $F_\gamma$ denotes its cumulative
distribution function, then $F_\gamma$ satisfies the equation
$$
 F(s) = \cases{
\sum_k p_k F^{*k} \left({\displaystyle{s \over 1-s}}
          \right), &if $s \in (0, 1)$; \cr
0, &if $s \le 0$; \cr
1, &if $s \ge 1$  . \cr }
\label e.cdf
$$
Moreover, the functional equation \ref e.cdf/ has exactly two solutions,
$F_\gamma$ and the Heaviside function $\I{\CO{0, \infty}}$ (see Lyons,
Pemantle and Peres (1997)).  Only a numerical computation of
$F_\gamma$ is known.  In all cases in which this
calculation was carried out, it appeared that the distribution was absolutely
continuous with respect to Lebesgue measure with bounded density.
Furthermore, unpublished work of the author and K.~Zumbrun has shown that
one can obtain {\it a posteriori\/} bounds on the modulus of continuity of
$F_\gamma$ from such computer calculations. Therefore, it seems reasonable
to conjecture, as in Lyons, Pemantle and Peres (1997),
that the distribution of $\gamma$ is always absolutely
continuous, at least when $p_k = 0$ for all sufficiently large $k$.

O.~H\"aggstr\"om (personal
communication, 1995) observed that one can write
$$
\gamma = [1, \gamma_1, 1, \gamma_2, 1, \gamma_3, \ldots]
:=
{1 \over \displaystyle {1 + {\mathstrut 1 \over \displaystyle \gamma_1 +
{\mathstrut 1 \over 1 + \cdots}}}}
\,,
$$
where $\gamma_i$ are i.i.d.\ In the special case $p_1 = p_2 = 1/2$,
the random variables $\gamma_i$ have the c.d.f.\ $(\I{\CO{0, \infty}} +
F_\gamma)/2$.  Note that $[1, x_1, 1, x_2, \ldots, 1, x_n]$ is monotonic
increasing in each $x_i$; this shows convergence of the infinite continued
fraction. Thus, in attempting to prove the above conjecture on absolute
continuity, 
we were led to consider the simpler random continued fractions
$$
[1, X_1, 1, X_2, 1, X_3, \ldots]\,,
$$
where $X_i$ are i.i.d.\ with the same distribution as
$$
X = \cases{0 &with probability 1/2,\cr
\alpha &with probability 1/2,\cr}
$$
where $\alpha \in (0, \infty)$ is a given real number. That is, instead of
trying to solve a fixed point problem, we are attempting to understand the
nature of the map taking a distribution of $X$ to the distribution of $[1,
X_1, 1, X_2, 1, X_3, \ldots]$ for a very simple $X$. Let $\mu_\alpha$ be
the distribution of this random continued fraction.

Our main result is the following. Define $\alpha_c$ to be the solution of
$\lambda_\alpha = {1 \over 2} \log 2$, where $\lambda_\alpha$ is defined in
\ref e.deflambda/ below.

\procl t.main
If $\alpha > 1/2$, then $\mu_\alpha$ is supported on a Cantor set of
Hausdorff dimension $< 1$. If $\alpha \le 1/2$, then the support of
$\mu_\alpha$ is an interval. If $\alpha > \alpha_c$, then
$\mu_\alpha$ is singular with respect to Lebesgue measure and is
concentrated on a set of Hausdorff dimension $< 1$; we have $\alpha_c \in
(0.2688, 0.2689)$. Whenever $\alpha > 0$, the measure $\mu_\alpha$ is
continuous.
\endprocl

The third sentence is the most novel one.

%As usual, by {\bf Hausdorff dimension} of a measure, $\mu$, we mean the
%smallest number $a$ such that there is a set $E$ of Hausdorff dimension $a$
%such that $\mu E^c=0$. (This is sometimes known as the upper Hausdorff
%dimension of $\mu$.)

Define
$$
T_\alpha (x) := {x+\alpha \over 1+x+\alpha}
\,.
$$
%and
%$$
%T := T_0\,.
%$$
Thus, $[1, x_1, 1, x_2, \ldots, 1, x_n] = T_{x_1} \circ T_{x_2} \circ
\cdots \circ T_{x_n} (0)$.
Note that $T_0$ is not a strict contraction.
The support of $\mu_\alpha$ is contained in $[0, M_\alpha]$ and
includes the endpoints, where $T_\alpha M_\alpha = M_\alpha > 0$.
If we solve this equation for $M_\alpha$,
we get $M_\alpha^2+\alpha M_\alpha-\alpha=0$, so that
$$
M_\alpha = {-\alpha+\sqrt{\alpha^2+4\alpha} \over 2}\,.
$$
In particular, $M_{1/2}=1/2$ and $T_\alpha 0 > T_0 M_\alpha$ iff $\alpha >
1/2$. Since $\mu_\alpha$ is supported on $T_0[0, M_\alpha] \cup T_\alpha [0,
M_\alpha]$, it follows by iteration that $\supp\, \mu_\alpha$ is a Cantor
set for $\alpha > 1/2$ and equals $[0, M_\alpha]$ for $\alpha \le 1/2$.
Thus,
\ref t.main/ asserts that in some interval of $\alpha$ where there
is overlap of $T_0[0, M_\alpha]$ and $T_\alpha [0, M_\alpha]$, the measure
$\mu_\alpha$ is singular. This contrasts with 
the case of two linear maps, $x \mapsto \alpha x$ and $x \mapsto
\alpha(1+x)$, where one has absolute continuity for almost
all $\alpha$ in the overlap region (Solomyak 1995, Peres and Solomyak
1996). 
However, based on some numerical evidence, we believe that there is also an
interval where $\mu_\alpha$ is absolutely continuous:

%Then $\mu_\alpha$ is the unique stationary probability measure for the
%Markov chain with transitions $s \mapsto T_X s$ (compare Diaconis and
%Freedman (1998)). 

\procl g.AC
For all $\alpha$ sufficiently small, $\mu_\alpha$ is absolutely continuous
with respect to Lebesgue measure.
\endprocl

After seeing a preprint of this work, Simon, Solomyak, and Urba\'nksi (1998)
made a great deal of progress on this conjecture. They proved that
for Lebesgue-a.e.\ $\alpha \in (0.215, \alpha_c)$, the measure
$\mu_\alpha$ is absolutely continuous. In particular, this shows that the
threshold $\alpha_c$ in \ref t.main/ is sharp.

Other work on random continued fractions includes Bernadac (1993, 1995),
Bhattacharya and Goswami (1998), Chamayou and Letac (1991), Chassaing,
Letac, and Mora (1984), Kaijser (1983), Letac (1986), Letac and Seshadri
(1995), Pitcher and Foster (1974), and Pincus (1983, 1985, 1994). In
particular, Pincus (1983) asks about singularity of distributions like
$\mu_\alpha$ when the images of two linear fractional maps overlap, which
is our main contribution.
Our analysis of the Hausdorff dimension is similar to that appearing on
pp.~166--168 of Bougerol and Lacroix (1985), where they use Lyapunov
exponent techniques to recover a result of Kinney and Pitcher (1965/1966).

Finally, returning to our original problem about Galton-Watson trees,
we mention another related
problem: Show that $\sum_{k\ge1} Z_k^{-1}$ has an absolutely continuous
distribution, where $Z_k$ is the size of the $k$th generation of a
Galton-Watson process. This arises as the resistance of the shorted
Galton-Watson tree.  In other words, if a graph is formed from a
genealogical Galton-Watson tree by identifying all vertices in level $k$
for each $k$ separately, then $\left(\sum_{k\ge1} Z_k^{-1}\right)^{-1}$ is
the effective conductance from the original progenitor to infinity. Again,
computation supports the conjecture that the distribution is absolutely
continuous with a bounded density, at least when the offspring distribution
is bounded.

In \ref s.proof/, we prove \ref t.main/. In \ref s.L^2/, we discuss the
possibility of $L^2$ or $L^\infty$ densities for $\alpha < \alpha_c$.

\bsection{Proof of \ref t.main/}{s.proof}

By Theorem 6.5 of Urba\'nski (1996), the Hausdorff dimension of $\supp\,
\mu_\alpha$ is less than 1 when there is no overlap, i.e., when
$\alpha > 1/2$.
(According to Remark 6.6 in that paper, we also have that the Hausdorff
dimension of $\supp\, \mu_\alpha$ equals its upper Minkowski dimension.)

Since $\mu_\alpha$ is a stationary measure for the Markov chain on $[0,
M_\alpha]$ that has transitions $s \mapsto T_X s$, it follows that
$\mu_\alpha$ is continuous: By stationarity, for any $s \in [0, M_\alpha]$,
we have
$$
\mu_\alpha(\{s\})
=
{1\over2} \mu_\alpha(\{T_0^{-1}s\}) + {1\over2} \mu_\alpha(\{T_\alpha^{-1}s\}) 
\,.
$$
In particular, if $s$ is such that
$\mu_\alpha(\{s\})$ is maximal,
then $\mu_\alpha(\{T_0^{-1}s\}) = \mu_\alpha(\{s\})$.
Likewise, $\mu_\alpha(\{T_0^{-n}s\}) = \mu_\alpha(\{s\})$.
Since $\mu_\alpha$ is finite, it follows that $\mu_\alpha(\{s\}) = 0$.

Define random variables $A_n, B_n, C_n, D_n$ by
$$
\mat(A_n, B_n, C_n, D_n) :=
\mat(1, X_1, 1, 1+X_1) \cdots \mat(1, X_n, 1, 1+X_n)\,.
$$
Since
$$
\mat(1, X, 1, 1+X)\vec(s, t) = \vec(s+Xt, s+t+Xt)
$$
and
$$
{s+Xt\over s+t+Xt} = {1\over \displaystyle 1 + {\mathstrut 1 \over \displaystyle
X+ {\mathstrut s \over t}}}\,,
$$
we have that
$$
{A_n s + B_n \over C_n s + D_n} =
[1, X_1, 1, X_2, \ldots, 1, X_n+s]\,,
$$
as is well known. We also write
$$
\mat(a_n, b_n, c_n, d_n) :=
\mat(1, x_1, 1, 1+x_1) \cdots \mat(1, x_n, 1, 1+x_n)
$$
for $x_i =0, \alpha$.

Let $Y$ be a random variable with distribution $\mu_\alpha$ independent of
all $X_i$. Then
$$
\mu_\alpha =
\L([1, X_1, \ldots, 1, X_n+Y])
$$
for any $n$,
where $\L(\cbuldot)$ is the law, i.e., the distribution, of a random
variable $\cbuldot$.
Now 
$$
\L([1, x_1, \ldots, 1, x_n+Y])
$$
is a probability measure supported on the interval
$$
\big[ [1, x_1, \ldots, 1, x_n],   [1, x_1, \ldots, 1, x_n+M_\alpha]\big]
=
\left[ {b_n \over d_n}, {a_n M_\alpha + b_n \over c_n M_\alpha + d_n}
\right]
\,.
$$
The length of this interval is
$M_\alpha/[d_n(c_n M_\alpha +d_n)]$ since
$$
a_n d_n - b_n c_n = \det\mat(a_n, b_n, c_n, d_n) = 1\,.
$$
Also,
$$
\mu_\alpha\Big(\big[ [1, x_1, \ldots, 1, x_n],   [1, x_1, \ldots, 1,
x_n+M_\alpha]\big]\Big) \ge 2^{-n}\,, \label e.msrbound
$$
with equality iff $\alpha \ge 1/2$. By Billingsley's Theorem (Falconer
(1997), p.~171), $\mu_\alpha$ is concentrated on a set of
Hausdorff dimension 
$$
\mu_\alpha\mathop{\hbox{-ess sup}}_x\ 
\liminf_{\epsilon \to0} {\log 1/\mu_\alpha\big((x-\epsilon,
x+\epsilon)\big) \over \log 1/(2\epsilon)}
\,,
\label e.dimmu
$$
which by \ref e.msrbound/ is at most
$$
\lim_{n \to\infty} {\log 2^n \over
\log \displaystyle {D_n(C_n M_\alpha +D_n) \over M_\alpha}}
=
{\log 2 \over 2\lambda_\alpha}\,,
\label e.dimbound
$$
where $\lambda_\alpha$ is the (top) Lyapunov exponent of the random matrix
$\mat(1, X, 1, 1+X)$; see \ref e.deflyap/ for the definition of
$\lambda_\alpha$ and Bougerol and Lacroix (1985), Cor.~VI.2.3, for this
property of $\lambda_\alpha$. [Simon, Solomyak and Urba\'nski (1998) have now
proved that the quantity in \ref e.dimmu/ is in fact equal to that in \ref
e.dimbound/ for Lebesgue-a.e.\ $\alpha \in (\alpha_c, 1/2)$.]
In particular, $\mu_\alpha$ is singular if
$\lambda_\alpha > {1\over2}\log 2$. Note
that $\lambda_\alpha$ is strictly increasing in $\alpha$ since $\mat(1, X,
1, 1+X)$ is stochastically increasing in $\alpha$. Thus, if $\alpha >
\alpha_c$, then $\mu_\alpha$ is singular and concentrated on a set of
Hausdorff dimension less than 1. It remains to estimate $\alpha_c$.

Parametrize projective one-space $\PP^1$ minus the horizontal direction by
the vectors $\R\times\{1\}$. Leaving out the horizontal is of no
consequence since $\mat(1, X, 1, 1+X)$ maps it to another direction
regardless of $X$.  The random map $v \mapsto \mat(1, X, 1, 1+X)v$ induces
a random map $(s, 1) \mapsto (T_X s, 1)$ on $\R\times\{1\}$ (thought of as
$\PP^1$).
The stationary measure for this Markov chain on
$\R\times\{1\}$ then becomes just $\mu_\alpha$ itself.  

Now $\lambda_\alpha$ is the expected change in the log norm on $\R^2$:
$$
\lambda_\alpha =
\int \E\log \big[\|\mat(1, X, 1, 1+X)(s, 1)\|/\|(s, 1)\|\big]
     \,d\mu_\alpha(s)\,.
\label e.deflyap
$$
If we use the norm
$$
\|(x, y)\|_\epsilon := \epsilon x + (1-\epsilon)y
$$
for any choice of $\epsilon \in (0, 1)$, then we get
$$
\lambda_\alpha
=
\int_{s=0}^{M_\alpha} \left\{
{1\over2} \log {\epsilon s + (1-\epsilon)(1+s) \over
\epsilon s + (1-\epsilon) } + 
{1\over2} \log {\epsilon (s+\alpha) + (1-\epsilon)(1+s+\alpha) \over
\epsilon s + (1-\epsilon) } \right\} \,d\mu_\alpha(s)\,.
\label e.formula
$$
Since this holds for all $\epsilon \in (0, 1)$, it also holds for $\epsilon
= 0$ by the Bounded Convergence Theorem:
$$
\lambda_\alpha =
{1\over2} \int_{s=0}^{M_\alpha} \log[(1+s)(1+s+\alpha)] \,d\mu_\alpha(s)
\,.
\label e.deflambda
$$

Define $K_\alpha$ to be the operator on c.d.f.'s on $[0, M_\alpha]$
corresponding to $\nu \mapsto ({1\over2}T_0+{1\over2}T_\alpha)\nu$. Thus,
$$
(K_\alpha F)(s) = {1\over2} F(T_0^{-1}s) + {1\over2}F(T_\alpha^{-1}s)
= {1\over2} F(T_0^{-1}s) + {1\over2}F(T_0^{-1}s-\alpha)\,.
$$
Set $F_0(s) := s/M_\alpha$ ($0 \le s \le M_\alpha$) and $F_{n+1} :=
K_\alpha F_n$. Then $F_n$ converges to the c.d.f.\ $F_{(\alpha)}$ of
$\mu_\alpha$ because $\left|\big[ [1, x_1, \ldots, 1, x_n],   [1, x_1,
\ldots, 1, x_n+M_\alpha]\big]\right| \le c_n^{-1} \le 1/n$.
Since $K_\alpha$ is a monotone operator (i.e., if $F \le G$,
then $K_\alpha F \le K_\alpha G$), we have that $F_{(\alpha)} \le F_{n+1} \le
F_n$ whenever $F_1 \le F_0$. If $\alpha < 1/2$, then
$F_1$ is linear from 0 to $T_\alpha 0$, from $T_\alpha 0$ to $T_0 M_\alpha$,
and from $T_0 M_\alpha$ to $M_\alpha$, with maximum slope in the middle
portion. If $\alpha \ge 1/2$, then $F_1$ is linear from 0 to $T_0 M_\alpha$,
constant from $T_0 M_\alpha$ to $T_\alpha 0$, and linear from $T_\alpha 0$ to
$M_\alpha$. Thus, in both cases,
it follows that $F_1 \le
F_0$ iff $F_1(T_0 M_\alpha) \le F_0(T_0 M_\alpha)$. Now
$$
F_0(T_0 M_\alpha) = {T_0 M_\alpha \over M_\alpha} = {1 \over 1+M_\alpha}
$$
and
$$
F_1(T_0 M_\alpha) = {1\over2} [F_0(M_\alpha) + F_0(M_\alpha-\alpha)]
=
1 - {\alpha \over 2M_\alpha}\,.
$$
Therefore, $F_1 \le F_0$ iff $\alpha \ge 1/6$.
(We remark that $M_{1/6} = 1/3$.)

Now if $F \le G$ and $h$ is an increasing function, then $\int h\,dF \ge
\int h\,dG$.
Since $s \mapsto 
\log[(1+s)(1+s+\alpha)]$ is increasing, this means that 
$$
\lambda_\alpha \ge {1\over2} \int \log[(1+s)(1+s+\alpha)]\,dF_n(s)
$$
for all $n$, giving us a lower bound for each $n$; in particular, as
computation shows, $\lambda_\alpha > {1\over2}\log 2$ for $\alpha =
0.2689$.
Therefore, $\alpha_c < 0.2689$.
%and hence for $\alpha \ge
%0.2689$. Thus, $\mu_\alpha$ is singular for such $\alpha$ and concentrated
%on a set of Hausdorff dimension $< 1$.
%\Qed

%\procl r.upperbd
We can get an upper bound on $\lambda_\alpha$ by a similar method. We have
that \ref e.formula/ holds for $\epsilon = 1$
by the Bounded and Monotone Convergence Theorems 
(use monotonicity in a neighborhood of $s=0$ and boundedness
elsewhere):
$$
\lambda_\alpha =
{1\over2} \int_{s=0}^{M_\alpha} \log\Big(1+{\alpha \over s}\Big)
       \,d\mu_\alpha(s)
\,.
$$
(Note that for best numerical estimates, though without error bounds, one
should choose a norm such that the expected log change is most nearly
constant, rather than $\epsilon =0, 1$ in \ref e.formula/. But the choices
$\epsilon = 0, 1$ give us lower and upper bounds when
$\alpha \ge 1/6$.)
\par
Now if $F \le G$ and $h$ is a decreasing function, then
$\int h\,dF \le \int h\,dG$. 
Since $s \mapsto \log(1+{\alpha/ s})$ is
decreasing, we get the upper bound
$$
\lambda_\alpha \le {1\over 2} \int \log\Big(1+{\alpha\over s}\Big)\,dF_n(s)
$$
for all $n$; in particular, $\lambda_\alpha < {1\over2}\log 2$ for $\alpha
\le 0.2688$.  
This completes the estimate of $\alpha_c$. \Qed
%\par
%We believe that $\mu_\alpha$ is singular for some $\alpha
%\le 0.2688$, but is absolutely continuous for sufficiently small $\alpha$.
%\endprocl

\bsection{$L^p$ Densities}{s.L^2}

We next show that $\mu_\alpha$ cannot be absolutely continuous with an
$L^2$ density for $\alpha > \sqrt 6/2-1 = 0.2247^+$. This result contrasts
with the case of two equally likely
linear maps (Solomyak 1995, Peres and Solomyak 1996),
where all known instances of absolutely continuous measures have an $L^2$
density.

\procl p.L^2
If $\mu_\alpha$ is absolutely continuous with an $L^2$ density, then
$\alpha \le \sqrt 6/2-1$.
\endprocl

\proof
Suppose that $\mu_\alpha$ has a density $f_\alpha$. Recall that $Y$ has
distribution $\mu_\alpha$. For $x_1, \ldots, x_n \in \{0, \alpha\}$,
the probability measure $\L([1, x_1, \ldots, 1, x_n+Y])$ has a density
$g_n$ (depending on $x_1, \ldots, x_n$) supported on the interval
$$
\left[ {b_n \over d_n}, {a_n M_\alpha + b_n \over c_n M_\alpha + d_n}
\right]
\,.
$$
Therefore,
$$
f_\alpha = \sum_{x_1, \ldots, x_n} 2^{-n} g_n\,.
$$
%If $f_\alpha \in L^2$, then
%$$
%\int_0^{M_\alpha} f_\alpha(s)^2\,ds < \infty\,.
%$$
Considering only the diagonal terms and using the Cauchy-Schwarz
inequality, we find that
\begineqalno
\int_0^{M_\alpha} f_\alpha(s)^2\,ds 
&\ge
\sum_{x_1, \ldots, x_n} 2^{-2n} \int g_n(s)^2 \,ds
\cr&\ge
\sum_{x_1, \ldots, x_n} 2^{-2n} \left|\left[ {b_n \over d_n}, {a_n M_\alpha
+ b_n \over c_n M_\alpha + d_n} \right]\right|^{-1}
\cr&=
\sum_{x_1, \ldots, x_n} 2^{-2n} {d_n(c_n M_\alpha + d_n) \over M_\alpha}
\cr&=
2^{-n} \Eleft{{D_n(C_n M_\alpha + D_n) \over M_\alpha}}
\cr&\ge
2^{-n}\E[D_n^2]/M_\alpha\,.
\cr
\endeqalno
Write
$$
\M_n := \mat(A_n, B_n, C_n, D_n)\,.
$$
Then 
$$
D_n^2 =
\pmatrix{0&0&0&1\cr} (\M_n\otimes \M_n) \pmatrix{0\cr0\cr0\cr1\cr}\,,
$$
whence
$$
\E[D_n^2]
=
\pmatrix{0&0&0&1\cr} \E[\M_n\otimes \M_n] \pmatrix{0\cr0\cr0\cr1\cr}\,.
$$
By independence, we have
$$
\Ebig{\M_n \otimes \M_n}
=
R^n\,,
$$
where
\begineqalno
R
&:=
\Eleft{\mat(1, X, 1, 1+X) \otimes \mat(1, X, 1, 1+X)}
\cr\noalign{\vskip5pt}&=
{1\over2} \tensormat(1, 0, 1, 0, 0, 1) +
{1\over2} \tensormat(1, \alpha, 1+\alpha, \alpha^2, \alpha + \alpha^2,
{(1+\alpha)^2})
\cr\noalign{\vskip8pt}&=
\tensormat(1, \alpha/2, 1+\alpha/2, \alpha^2/2, \alpha/2+\alpha^2/2,
1+\alpha+\alpha^2/2)\,.
\cr
\endeqalno
The characteristic polynomial of $R$ is
$$
t\mapsto
t^4-(4+2\alpha +\alpha^2/2)t^3 + (6+4\alpha+\alpha^2/2)t^2-(4+2\alpha)t+1
\,,
$$
which has its largest root $> 2$ iff $\alpha > \sqrt 6/2-1$.  Clearly, the
Perron-Frobenius eigenvector of $R$ has all its coordinates strictly
positive, so that 
$$
2^{-n}\E[D_n^2]
=
2^{-n} \pmatrix{0&0&0&1\cr} R^n \pmatrix{0\cr0\cr0\cr1\cr}
\to\infty
$$
as $n \to\infty$ if $\alpha > \sqrt 6/2-1$.  Hence $\mu_\alpha$ cannot
have a density in $L^2$ for $\alpha > \sqrt 6/2-1$.
\Qed

Similar methods show that if $\alpha > 3\sqrt 2 - 4 = 0.24^+$,
then $\mu_\alpha$ cannot have a density in $L^{3/2}$.
In the other direction, we can extend the method to show the following
proposition:

\procl p.L^p If $\mu_\alpha$ is absolutely continuous with a density in
$L^p$ for every $p < \infty$, then $\alpha \le (3 \sqrt 2 - 4)/2 =
0.1213^+$. In particular, this is the case if $\mu_\alpha$ has a density
in $L^\infty$.
\endprocl

\proof
Let $p$ be a positive integer and
$$
R_p(\alpha) :=
\Eleft{\M_1^{\otimes 2(p-1)}}
\,,
$$
where $\M_1^{\otimes r}$ denotes the $r$-th tensor power of
$\M_1 := \mat(1, X, 1, 1+X)$. Write $r := 2(p-1)$.
The method of proof of \ref p.L^2/ shows that
$\mu_\alpha$ cannot have a density in $L^p$ for $\alpha > \alpha_p$, where
$\alpha_p$ is the value of $\alpha$ such that the largest eigenvalue
$\gamma_p$ of
$R_p$ is equal to $2^{(p-1)/2}$.
Let the eigenvectors of $\M_1$ be $\u_1$ and $\u_2$. Then the
eigenvectors of $\M_1^{\otimes r}$ are $\u_{i_1} \otimes \cdots \otimes
\u_{i_r}$, where all $i_j$ are either 1 or 2. Thus, the largest eigenvalue
of $\M_1^{\otimes r}$ is $\delta^r$, where $\delta$ is the largest
eigenvalue of $\M_1$. 
It follows that $\gamma_p^{1/r}$ tends to the largest eigenvalue of
$\mat(1, \alpha, 1, 1+\alpha)$. It is then an easy matter to show that
$\alpha_p \to (3 \sqrt 2 - 4)/2$.
\Qed

\procl r.bdddensity Numerical evidence suggests that in fact $\mu_\alpha$
may not have a bounded density until $\alpha$ is less than about 0.05 (very
roughly). Furthermore, the density seems to gain increasing smoothness the
smaller $\alpha$ becomes.
\endprocl

%it suffices to show that $2^{-n} \Ebig{\| \M_n \otimes \M_n
%\|_1} \to \infty$, where the $\|\cbuldot\|_1$-norm of a matrix is the sum
%of the absolute values of its entries, since $D_n(C_n M_\alpha + D_n) /
%M_\alpha \ge D_n C_n \ge \| \M_n \otimes \M_n\|_1$ [NO -- it's $\le$!]. 

%Now for independent matrices $\M$, $\M'$ with nonnegative coefficients,
%we have 
%$$
%\Ebig{\|\M \M' \|_1} = \|(\E\M) (\E\M')\|_1\,.
%$$
%Therefore,
%$$
%\Ebig{\| \M_n \otimes \M_n \|_1}
%=
%\left\|\left(\Eleft{\mat(1, X, 1, 1+X) \otimes \mat(1, X, 1,
%1+X)}\right)^n\right\|_1
%\,.
%$$
%We have
%\begineqalno
%\Eleft{\mat(1, X, 1, 1+X) \otimes \mat(1, X, 1, 1+X)}
%&=
%{1\over2} \tensormat(1, 0, 1, 0, 0, 1) +
%{1\over2} \tensormat(1, \alpha, 1+\alpha, \alpha^2, \alpha + \alpha^2,
%{(1+\alpha)^2})
%\cr&=
%\tensormat(1, \alpha/2, 1+\alpha/2, \alpha^2/2, \alpha/2+\alpha^2/2,
%1+\alpha+\alpha^2/2)\,.
%\cr
%\endeqalno
%The characteristic polynomial of this matrix is
%$$
%t\mapsto
%t^4-(4+2\alpha +\alpha^2/2)t^3 + (6+4\alpha+\alpha^2/2)t^2-(4+2\alpha)t+1
%\,,
%$$
%which has its largest root $> 2$ iff $\alpha > \sqrt 6/2-1$.

%$$
%\lim_{n \to\infty} \L([1, X_n, \ldots, 1, X_1+U])
%=
%\lim_{n \to\infty} ({1\over2}T_0 +{1\over2}T_\alpha)^n \L(U)\,.
%$$
%
%$$
%For any choice of $\epsilon_i = 0, \alpha$, we have
%$$
%T_{\epsilon_1} T_{\epsilon_2} \cdots T_{\epsilon_n} \L(U)
%=
%\U[T_{\epsilon_1} \cdots T_{\epsilon_n} 0, T_{\epsilon_1} \cdots
%T_{\epsilon_n} M\alpha]
%=
%\U\big[ [1, \epsilon_1, 1,\epsilon_2, \ldots, 1, \epsilon_n],
%[1, \epsilon_1, 1,\epsilon_2, \ldots, 1, \epsilon_n + M_\alpha]\big]\,.
%$$

\medbreak
\noindent {\bf Acknowledgement.}\enspace I thank Grahame Bennett,
Rick Kenyon, Yuval Peres, and Kevin Zumbrun for useful conversations. I
also thank the referee for suggesting the investigation of $L^\infty$
densities. I am grateful to the Universit\'e de Lyon for its hospitality,
where most of this work was carried out. 

\beginreferences

Bernadac, E. (1993)
Fractions continues al\'eatoires sur un cone sym\'etrique, {\it
C. R. Acad. Sci. Paris}, I, {\bf 316}, 859--864.

Bernadac, E. (1995) Random continued fractions and inverse Gaussian
distribution on a symmetric cone, {\it J. Theoret. Probab.} {\bf 8},
221--259.

Bhattacharya, R. \and Goswami, A. (1998) A class of random continued
fractions with singular equilibria, {\it preprint}.

Bougerol, P. \and Lacroix, J. (1985) {\it Products of Random Matrices with
Applications to Schr\"odinger Operators.} Birkha\"user, Boston.

Chamayou, J.-F. \and Letac, G. (1991) Explicit stationary distributions for
compositions of random functions and products of random matrices, {\it J.
Theoret. Probab.} {\bf 4}, 3--36.

Chassaing, P., Letac, G. \and Mora, M. (1984) Brocot sequences and random
walks in ${\rm SL}(2,{\bf R})$, in {\it Probability Measures on Groups, VII
(Oberwolfach, 1983)}, 36--48, Lecture Notes in Math., {\bf 1064}, Springer,
Berlin.

Doyle, P. G. \and Snell, J. L., (1984) {\it Random Walks and
Electric Networks.}  Mathematical Assoc. of America, Washington, DC.

Falconer, K. (1997) {\it Techniques in Fractal Geometry.} Wiley,
Chichester.

Kaijser, T. (1983) A note on random continued fractions, in {\it
Probability and Mathematical Statistics}, Uppsala Univ., Uppsala, pp.
74--84.

Kinney, J. R. \and Pitcher, T. S. (1965/1966) The dimension of some sets
defined in terms of $f$-expansions, {\it Z. Wahrscheinlichkeitstheorie und
Verw. Gebiete} {\bf 4}, 293--315.

Letac, G. (1986) A contraction principle for certain Markov chains and its
applications, in {\it Random Matrices and Their Applications (Brunswick,
Maine, 1984)}, Contemp. Math., {\bf 50}, pp.  263--273. Amer. Math. Soc.,
Providence, R.I.

Letac, G. \and Seshadri, V. (1995) A random continued fraction in ${\bf
R}\sp {d+1}$ with an inverse Gaussian distribution, {\it Bernoulli\/} {\bf
1}, 381--393.

Lyons, R., Pemantle, R. \and Peres, Y. (1995) Ergodic theory on
Galton-Watson trees: speed of random walk and dimension of harmonic
measure, {\it Ergodic Theory Dynamical Systems} {\bf 15}, 593--619.

Lyons, R., Pemantle, R. \and Peres, Y. (1997) Unsolved problems concerning
random walks on trees, in {\it Classical and Modern Branching Processes}, K.
Athreya and P. Jagers (editors), pp.~223--238. Springer, New York.

Lyons, R. \and Peres, Y. (1998) {\it Probability on Trees and Networks}.
Cambridge University Press, in preparation. Current version available at
\hfill\break {\tt http://php.indiana.edu/\~{}rdlyons/}.

Peres, Y. \and Solomyak, B. (1996) Absolute continuity of Bernoulli
convolutions, a simple proof, {\it Math.  Res. Lett.} {\bf 3}, 231--239.

Pincus, S. (1983) A class of Bernoulli random matrices with continuous
singular stationary measures, {\it Ann. Probab.} {\bf 11}, 931--938.

Pincus, S. (1985) Strong laws of large numbers for products of random
matrices, {\it Trans. Amer.  Math. Soc.} {\bf 287}, 65--89.

Pincus, S. (1994) Singular stationary measures are not always fractal, {\it
J.  Theoret. Probab.} {\bf 7}, 199--208.

Pitcher, T. \and Foster, S. (1974) Convergence properties of random
$T$-fractions, {\it Z. Wahrscheinlichkeitstheorie und Verw. Gebiete} {\bf
29}, 323--330.

Simon, K., Solomyak, B. \and Urba\'nski, M. (1998) Parabolic iterated
function systems with overlaps II: invariant measures, {\it preprint}.

Solomyak, B. (1995) On the random series $\sum \pm\lambda^n$ (an Erd\"os
problem), {\it Ann. of Math.} (2) {\bf 142}, 611--625.

Urba\'nski, M. (1996) Parabolic Cantor sets, {\it Fund. Math.} {\bf 151},
241--277.

\endreferences

\filbreak
\begingroup
\eightpoint\sc
\parindent=0pt\baselineskip=10pt

\def\emailwww#1#2{\par\qquad {\tt #1}\par\qquad {\tt #2}\smallskip}

Department of Mathematics,
Indiana University,
Bloomington, IN 47405-5701, USA
\emailwww{rdlyons@indiana.edu}
{http://php.indiana.edu/\~{}rdlyons/}

\endgroup

\bye